\documentclass{amsart}

 \usepackage[top=1.65in, bottom=1.65in, right=1.3in, left=1.3in]{geometry} 

\usepackage[francais,english]{babel}

\usepackage{tikz}
\usetikzlibrary{matrix,arrows,calc,shapes,trees,positioning}

\usepackage{amsfonts}
\usepackage{amssymb}
\usepackage{amsmath}
\usepackage{fancyhdr}
\usepackage{eufrak}
\usepackage{latexsym}
\usepackage{amsbsy}
\usepackage[all]{xy}
\usepackage[latin1]{inputenc}
\usepackage[T1]{fontenc}        
\usepackage{amsmath}

\newtheorem{defn}{D\'efinition}[section]
\newtheorem{thm}[defn]{Th\'eor\`eme}
\newtheorem{prop}[defn]{Proposition}
\newtheorem{lem}[defn]{Lemme}

\newtheorem{rem}[defn]{Remarque}

\newtheorem{note}[defn]{Notations}

\renewcommand{\proof}{\vskip 2mm \noindent {\textsc{D\'emonstration: }}\rm}

\newcommand{\fin}{\hfill{\Large$\Box$}\\}

\newcommand{\la}{\lambda}

\newcommand{\C}{\mathbb {C}}

\newcommand{\N}{\mathbb {N}}
\newcommand{\Z}{\mathbb {Z}}

\newcommand{\Pj}{\mathbb {P}}

\def\com{\ar@{}[rd]|{\circlearrowleft}}

\DeclareMathOperator{\Bif}{Bif}

\title {Perturbations d'exemples de Latt\`es et dimension de Hausdorff du lieu de bifurcation}

\thanks{This research was partially supported by the ANR project LAMBDA, ANR-13-BS01-0002}

\begin{author}[F.~Berteloot]{Fran\c cois Berteloot}
\address{ 
 Universit\'e de Toulouse - IMT\\
 UMR CNRS 5219\\
 31062 Toulouse Cedex \\
  France }
  \email{francois.berteloot$@$math.univ-toulouse.fr}
\end{author}

\begin{author}[F.~Bianchi]{Fabrizio Bianchi}
\address{ 
 Imperial College\\
 South Kensington Campus \\
  London SW7 2AZ\\
  UK }
  \email{f.bianchi$@$imperial.ac.uk}
\end{author}

\begin{document}

\maketitle

\begin{abstract}
We give an estimate for the Hausdorff dimension of the bifurcation locus of a family of endomorphisms of $\Pj^k (\C)$. This dimension is maximal near isolated Latt\`es examples.
\end{abstract}

\section{Introduction}

Un \emph{exemple de Latt\`es} de degr\'e alg\'ebrique $d$ est un endomorphisme holomorphe $f$ de $\Pj^k$
faisant commuter un diagramme
\begin{center}
\begin{tikzpicture}
\matrix (m) [matrix of math nodes, row sep=3.5em,
column sep=3.5em, text height=1.5ex, text depth=0.5ex]
{ \C^k& \C^k \\
 \Pj^k & \Pj^k \\};
	\path[->,font=\small]
		(m-1-1) edge node[left] {$\sigma$} (m-2-1)
		(m-1-2) edge node[auto] {$\sigma$} (m-2-2)
		(m-1-1) edge node[auto] {$ D $} (m-1-2)
		(m-2-1) edge node[auto] {$ f $} (m-2-2);
\end{tikzpicture}
\end{center}
\noindent o\`u $D$ est une application affine de partie lin\'eaire $\sqrt{d}U$ ($U$ unitaire) et $\sigma$ est un rev\^{e}tement ramifi\'e sur les fibres duquel un groupe cristallographique complexe agit transitivement. Ces exemples ont \'et\'es d\'ecouverts en 1918 par S. Latt\`es   pour la dimension $k=1$ \cite{Lat} mais ils existent en toute dimension et tout degr\'e. L'article  de J. Milnor \cite{Mil} d\'ecrit ces objets d'un point de vue contemporain.
 La caract\'erisation suivante, due \`a A. Zdunik \cite{Zd} en dimension $k=1$  et F. Berteloot, C. Dupont et J.J. Loeb \cite{BL}, \cite{BDu},\cite{Du} en toute dimension,  permet d'amorcer l'\'etude des bifurcations engendr\'ees par l'extr\`eme rigidit\'e de ces endomorphismes.

\begin{thm}\label{THBDL}Soit $f$ un endomorphisme holomorphe de degr\'e   alg\'ebrique
$d$ sur $\Pj^k$ et $\mu_f$ sa mesure d'entropie maximale. Les assertions suivantes sont \'equivalentes :
\begin{itemize}
\item[1)] $\mu_f$ est absolument continue par rapport \`a la mesure de Lebesgue,
\item[2)] les exposants de Lyapounov de $(f,\mu_f)$ sont  minimaux \'egaux \`a 
$\ln \sqrt{d}$,
\item[3)] $f$ est un exemple de Latt\`es.
\end{itemize}
\end{thm} 

Une \emph{famille holomorphe} d'endomorphismes de $\Pj^k$, param\'etr\'ee par une vari\'et\'e complexe $M$, est une application  holomorphe
$F: M\times \Pj^k\to M\times \Pj^k$ de la forme $(\la,f_\la(z))$ et telle que 
le degr\'e alg\'ebrique des endomorphismes $f_\la$ soit \'egal \`a $d$ pour tout 
$\la\in M$. En dimension $k=1$, la nature des bifurcations au sein de telles familles est bien comprise depuis les travaux fondateurs de R. Ma\~{n}\'e, P. Sad et D. Sullivan \cite{MSS} et, ind\'ependemment, M. Lyubich \cite{L1},\cite{L2}. De plus, comme l'a montr\'e L. DeMarco \cite{dM}, le lieu de bifurcation ${\Bif(F)} $ co\"{i}ncide avec
le support du $(1,1)$-courant positif ferm\'e $dd^c L(\la)$ o\`u $L(\la)$
d\'esigne l'exposant de Lyapounov de $f_\la$ par rapport \`a sa mesure d'entropie
maximale. Dans un travail en collaboration avec  C. Dupont, les auteurs  ont \'etendu cette th\'eorie en dimension arbitraire. Dans ce contexte la fonction $L(\la)$ d\' esigne la somme des exposants de Lyapounov de $f_\la$ par rapport \`a sa mesure d'entropie maximale $\mu_\la$, l'ensemble de Julia $J(\la)$ de $f_\la$
est par d\'efinition le support de $\mu_\la$ et un $J$-cycle est un cycle contenu dans $J(\la)$. Le principal r\'esultat est le suivant, nous renvoyons \`a  \cite{BBD} pour un \'enonc\'e plus complet.

 \begin{thm}\label{THBBD}  
Soit $F:M\times \Pj^{k}\to M\times \Pj^{k}$  
une famille holomorphe d'endomorphismes o\`u $M$ est un ouvert simplement connexe de l'espace des endomorphismes
de $\Pj^k$ de degr\'e $d \geq 2$. 
Alors les assertions suivantes sont  \'equivalentes :
\begin{enumerate}
\item[1)] les $J$-cycles r\'epulsifs de $f_\la$ bougent holomorphiquement sur $M$,
\item[2)]  la fonction $L$ est pluriharmonique sur $M$,
\item[3)]  $J(\la)$ bouge holomorphiquement sur $M$.
\end{enumerate}
\end{thm} 

On d\'efinit alors le lieu de bifurcation  ${\Bif(F)}$  d'une telle famille comme
le support du $(1,1)$-courant positif ferm\'e $dd^c L(\la)$. Il est remarquable qu'en dimension $k\ge 2$ (et contrairement au cas de la dimension $k=1$) le lieu de bifurcation puisse \^{e}tre d'int\'erieur non vide.
Ce ph\'enom\`ene a \'et\'e r\'ecemment mis en \'evidence par le deuxi\`eme auteur et J. Taflin \cite{BiTa} et R. Dujardin
\cite{Duj}.\\

Lorsque $f_{\la_0}$ est un exemple de Latt\`es isol\'e on voit,  gr\^{a}ce \` a la caract\'erisation par la minimalit\'e de la somme des exposants de Lyapounov, que $\la_0$  est un param\`etre de bifurcation. Nous pensons
que ces endomorphismes sont en fait le foyer de bifurcations maximales qu'il conviendrait d'\'etudier.

 Illustrons ceci en rappelant ce qui est connu en dimension $k=1$. Comme l'ont remarqu\'e G. Bassanelli et F. Berteloot, l'argument de minimalit\'e de l'exposant montre aussi que $\la_0$ est dans le support des puissances ext\'erieures de $dd^c L$ ce qui permet de pr\'eciser le nombre de cycles r\'epulsifs bifurquants simultan\'ement. 
Par exemple, dans l'espace des fractions rationnelles de degr\'e $d$, on peut montrer qu'un exemple de Latt\`es non flexible est \`a la fois accumul\'e par des fractions rationnelles hyperboliques et par des fractions rationnelles poss\'edant $2d-2$ cycles neutres distincts
(\cite{BB1}, voir aussi \cite{Be} subsection 6.2.3). Ceci a \'et\'e \'etendu aux exemples de Latt\`es flexibles par X. Buff et T. Gauthier \cite{BG}. 
Un r\'esultat important de M. Shishikura stipule qu'une fraction de degr\'e $d$ poss\`ede au plus $2d-2$ cycles non r\'epulsifs et que cette borne est r\'ealis\'ee \cite{Shi}. Les bifurcations  de ce type  sont donc, en un certain sens, maximales. Soulignons que T. Gauthier \cite{Ga} a montr\'e que la dimension de Hausdorff du lieu des bifurcations maximales est elle-m\^{e}me maximale.

Nous initions dans cet article l'\'etude des bifurcations engendr\'ees par un exemple  de Latt\`es en dimension quelconque. Nous montrons que le lieu de ces bifurcations est de dimension de Hausdorff maximale dans "toutes les directions".
Notre principal r\'esultat est le suivant.

\begin{thm}\label{thP}Soit $F:D\times \Pj^k \to D\times \Pj^k$ une famille holomorphe d'endomorphismes de degr\'e $d\ge 2$ param\'etr\'ee par le disque unit\'e $D$ de $\C$.
On suppose que $f_0$ est un exemple de Latt\`es et que $0$ est accumul\'e par des param\`etres $\la\in D$ tel que $f_\la$ n'est pas un exemple 
de Latt\`es. Alors $\textrm{dim}_H\left( {\Bif(F)} \right) = 2$.
\end{thm}

Terminons cette introduction en pr\'ecisant quelques notations.

\begin{note} $D$ est le disque unit\'e de $\C$ et  $D_r$ d\'esigne le disque $rD$.\\
Pour tout sous-ensemble $E$ d'un produit $D\times B$ et tout $\la\in D$, on note $(E)_\la$ la tranche 
$E\cap(\{\la\}\times B)$. On note $\pi_D$ la projection canonique sur $D$.\\
 On note $\Gamma_\gamma:=\{(\la,\gamma(\la))\;:\; \la\in D\}$ le graphe d'une application $\gamma : D\to B$.\\
${\mathcal O}\left(D,B\right)$ est l'espace des applications holomorphes de $D$ dans $B$.\\
${\mathcal H}_d\left(\Pj^k\right)$ est l'espace des endomorphismes holomorphes de degr\'e alg\'ebrique $d$ sur $\Pj^k$.
\end{note}

\section{Laminations engendr\'ees par des contractions}\label{Sec1}
Dans toute cette section, $B$ d\'esignera une boule de $\C^k$ pour la norme hermitienne standard. Rappelons que 
 le disque unit\'e de $\C$ est not\'e $D$.
Nous allons construire une lamination dans $D\times B$  par une famille  $\mathcal G$ de graphes
holomorphes au-dessus de $D$ dont les
 tranches  sont  des ensembles de Cantor. La dimension de Hausdorff des tranches sera  minor\'ee en utilisant les travaux de Pesin et Weiss  \cite{PW}
 sur les constructions g\'eom\'etriques du type Moran (voir aussi \cite[Chapter 5]{Pe})
.\\

 L'ensemble $\mathcal G$ est obtenu par un proc\'ed\'e usuel \`a partir d'une collection $G_1,\cdots,G_m$ de 
contractions holomorphes de la forme $G_j(\la,z)=\left(\la,G_{j,\la}(z)\right)$, d\'efinies sur un voisinage $\widetilde{D}\times \widetilde{B}$ de $\overline{D}\times\overline{B}$ et v\'erifiant les propri\'et\'es suivantes pour des constantes 
$ 0<a\le A$

\begin{itemize}
\item[1)] $G_j \left(\overline{D}\times\overline{B}\right) \subset \overline{D}\times{B},\; \forall j\in\{1,\cdots,m\}$
\item[2)] $\textrm{dist}\left(G_j \left(\overline{D}\times\overline{B}\right), G_k \left(\overline{D}\times\overline{B}\right)\right) >0,\; \forall j\ne k\in\{1,\cdots,m\}$
\item[3)]  $e^{-A} \Vert z-z'\Vert \le \Vert G_j(\la,z)-G_j(\la,z')\Vert \le e^{-a} \Vert z-z'\Vert,\;
\forall \la \in \overline{D}, \forall z,z'\in \overline{B}, \forall j\in\{1,\cdots,m\}$.
\end{itemize}

Pour tout $\omega:=(\omega_k)_{k\ge 0} \in \{1,\dots,m\}^{\N}=:\Sigma_m^+$, on pose $G_{\omega_0\cdots \omega_p}:=G_{\omega_0}\circ \cdots \circ G_{ \omega_p}$ puis
\begin{center}
$T_{\omega_0\cdots \omega_p}:=G_{\omega_0\cdots \omega_p}\left( \overline{D}\times\overline{B}\right)$
\end{center}
\begin{center}
$\Gamma_\omega := \cap_{p\ge 0} T_{\omega_0\cdots \omega_p}$.
\end{center}
On observera que les applications $G_{\omega_0\cdots \omega_p}$ sont de la forme $(\la,z)\mapsto (\la,G_{ \omega_0\cdots \omega_p,\la}(z))$ et que
\begin{eqnarray}\label{a}
\Vert G_{\omega_0\cdots \omega_p}(\la,z) -G_{\omega_0\cdots \omega_p}(\la,z')\Vert\le e^{-(p+1)a} \Vert z-z'\Vert,
\; \forall \la \in \overline {D}, \forall z,z'\in \overline{B}.
\end{eqnarray} 
En particulier, $\left(T_{\omega_0\cdots \omega_p}\right)_\la$ est une suite d\'ecroissante de compacts dont le diam\`etre tend vers $0$ lorsque $p$ tend vers $+\infty$ et il existe donc un unique point $\omega(\la)\in B$ tel que
\begin{center}
$\left(\cap_{p\ge 0} T_{\omega_0\cdots \omega_p}\right)_\la = \cap_{p\ge 0}\left(T_{\omega_0\cdots \omega_p}\right)_\la=\{(\la,\omega(\la))\}.$
\end{center}
Autrement dit, $\Gamma_\omega$ est le graphe d'une application  d\'efinie sur $\overline{D}$ et \`a valeurs dans $B$ que l'on note aussi $\omega$. On notera $\mathcal{G}$ la r\'eunion des graphes ainsi obtenus
\begin{center}
${\mathcal G}:= \bigcup_{\omega\in \Sigma_m^+}\Gamma_\omega$.
\end{center}

\begin{prop}\label{PropCantGr}
L'ensemble $\mathcal G$ est  constitu\'e de graphes deux \`a deux disjoints, continus sur $\overline D$ et holomorphes sur $D$. 
Pour tout $\la_0 \in D$, la dimension de Hausdorff de $\left({\mathcal G}\right)_{\la_0}$ est minor\'ee par $\frac{\ln m}{A}$ et l'application $H_{\la_0} : {\mathcal G}\to \left({\mathcal G}\right)_{\la_0}$ d\'efinie par $H_{\la_0}(\la,\omega(\la))=(\la_0,\omega(\la_0))$ est $\frac{a}{A}$-H\"{o}lder sur ${\mathcal G}\cap\left(D_r\times B\right)$ pour tout 
$0<r<1$.
\end{prop}

\proof Fixons $z_0\in B$. L'in\'egalit\'e (\ref{a}) appliqu\'ee \`a $z'=G_{\omega_{p+1}\cdots\omega_{p+q}}(\la,z_0)$ et $z=z_0$   donne
\begin{eqnarray}\label{b}
\Vert G_{\omega_0\cdots \omega_p}(\la,z_0) -G_{\omega_0\cdots \omega_{p+q}}(\la,z_0)\Vert\le e^{-(p+1)a} \textrm{diam}(B), \; \forall \la\in \overline{D}.
\end{eqnarray} 
Les propri\'et\'es de r\'egularit\'e de $\omega$ r\'esultent alors de la convergence uniforme sur $\overline D$ de $G_{\omega_0\cdots \omega_p}(\la,z_0)$ vers $(\la,\omega(\la))$.\\

 Nous allons maintenant montrer  que pour tout $0<r<1$ il existe une constante $C_r>0$ telle que si $(\omega_0,\cdots,\omega_p)\ne (\omega'_0,\cdots,\omega'_p) $ alors
\begin{eqnarray}\label{c}
\Vert G_{\omega_0\cdots \omega_p}(\la,z) -G_{\omega'_0 \cdots \omega'_{p}}(\la',z')\Vert\ge C_re^{-pA}, \; \forall (\la,z),(\la',z')\in D_r\times B.
\end{eqnarray} 
Puisque $G_{\omega_0\cdots \omega_p}(\overline{D}\times\overline{B})\subset {\overline D}\times B$, on d\'eduit des in\'egalit\'es de Cauchy qu'il existe une constante $K_r\ge 1$ ind\'ependante de $(\omega_0,\cdots,\omega_p)$ telle que 
\begin{eqnarray}\label{d}
\Vert G_{\omega_0\cdots \omega_p}(\la,z) -G_{\omega_0\cdots \omega_p} (\la',z) \Vert\le K_r \vert \la-\la'\vert, \; \forall (\la,z),(\la',z)\in D_r\times \overline{B}.
\end{eqnarray}
Soit $0<d:=\min_{1\le j\ne k\le m} \textrm{dist}\left(G_j \left(\overline{D}\times\overline{B}\right), G_k \left(\overline{D}\times\overline{B}\right)\right)$. Soit $(\la,z),(\la,z') \in D_r\times\overline{B}$, comme $\Vert G_{\omega_j}(\la,z)-G_{\omega'_j}(\la,z')\Vert \ge e^{-A}\Vert z-z'\Vert$ si $\omega_j = \omega'_j$
et $\Vert G_{\omega_j}(\la,z)-G_{\omega'_j}(\la,z')\Vert \ge d$ sinon, on voit que
\begin{eqnarray}\label{e}
\Vert G_{\omega_0\cdots \omega_p}(\la,z) -G_{\omega'_0 \cdots \omega'_{p}}(\la,z')\Vert\ge e^{-pA}d, \;\; \forall (\la,z),(\la,z')\in D_r\times B.
\end{eqnarray}
Supposons que  $\vert \la-\la'\vert \le \frac{d}{2K_r} e^{-pA}$,  on d\'eduit alors de (\ref{d}) et (\ref{e}) que la minoration   
$\Vert G_{\omega_0\cdots \omega_p}(\la,z) -G_{\omega'_0 \cdots \omega'_{p}}(\la',z')\Vert \ge e^{-pA} d -K_r \vert \la-\la'\vert \ge \frac{d}{2} e^{-pA}
\ge \frac{d}{2K_r} e^{-pA}$ a lieu pour tout $z,z'\in B$. Comme 
$\Vert G_{\omega_0\cdots \omega_p}(\la,z) -G_{\omega'_0 \cdots \omega'_{p}}(\la',z')\Vert \ge \vert \la-\la'\vert$, la m\^{e}me minoration reste vraie lorsque $\vert \la-\la'\vert \ge \frac{d}{2K_r} e^{-pA}$. Ceci justifie la minoration (\ref{c})
avec $C_r:=\frac{d}{2K_r}$.\\

Nous pouvons maintenant terminer la preuve de la proposition. Soit $\la_0\in D$, l'in\'egalit\'e (\ref{c}) montre en particulier qu'il existe une constante $C>0$ que 
\begin{eqnarray*}
\textrm{dist}\left((T_{\omega_0\cdots\omega_p})_{\la_0}, (T_{\omega'_0\cdots\omega'_p})_{\la_0}\right)\ge C e^{-pA}\;\textrm{ si }\;(\omega_0,\cdots,\omega_p)\ne (\omega'_0,\cdots,\omega'_p).
\end{eqnarray*}
 Un th\'eor\`eme d\^{u} \`a Pesin et Weiss  (voir \cite{PW}, Proposition 5) stipule   que dans ces conditions $\textrm{dim}_H \left(\left({\mathcal G}\right)_{\la_0}\right) \ge \frac{\ln m}{A}$.

Il nous reste \`a \'etudier la r\'egularit\'e de l'application $H_{\la_0}$. Soient $\omega \ne \omega' \in \Sigma_m^+$ et $p:=\min \{j\;:\; \omega_j\ne \omega'_j\}$. Fixons $z_0\in B$ et, pour $q>p$, posons
$(\la,z):=G_{\omega_{p+1}\cdots\omega_q}(\la,z_0)$, $(\la',z'):=G_{\omega'_{p+1}\cdots\omega'_q}(\la',z_0)$. D'apr\`es (\ref{c}) il vient
$\Vert G_{\omega_0\cdots \omega_q}(\la,z_0) -G_{\omega'_0 \cdots \omega'_{q}}(\la',z_0)\Vert=
\Vert G_{\omega_0\cdots \omega_p}(\la,z) -G_{\omega'_0 \cdots \omega'_{p}}(\la',z')\Vert \ge C_re^{-pA}$ pour tout $\la,\la' \in D_r$ d'o\`u, en faisant tendre 
$q$ vers $+\infty$,
\begin{eqnarray}\label{f}
\Vert (\la,\omega(\la)) - (\la',\omega'(\la'))\Vert \ge C_r e^{-pA},\;\;\forall\la,\la' \in D_r.
\end{eqnarray}
Par ailleurs, puisque $(\la_0,\omega(\la_0))$ et $(\la_0,\omega'(\la_0))$ sont tous deux dans $G_{\omega_0\cdots \omega_{p-1}}(D\times B)$, on a 
\begin{eqnarray}\label{g}
\Vert \omega(\la_0) - \omega'(\la_0)\Vert \le  e^{-pa} (\textrm{diam}\;B).
\end{eqnarray}
On tire  de (\ref{f}) et (\ref{g}) que $\Vert \omega(\la_0)-\omega'(\la_0)\Vert \le C_r^{-\frac{a}{A}} (\textrm{diam}\;B)   \Vert (\la,\omega(\la)) - (\la',\omega'(\la'))\Vert^{\frac{a}{A}}$ pour $\la,\la'\in D_r$ ce qui signifie que l'application $H_{\la_0}$ est $\frac{a}{A}$-H\"{o}lder sur ${\mathcal G}\cap\left(D_r\times B\right)$.\fin

Supposons maintenant que $D\times B$ contient une hypersurface irreductible $Z$ qui n'est pas r\'eduite \`a une fibre de $\pi_D$ et pour laquelle $\pi_D(Z) \subset D_{r_0}$ o\`u $0<r_0<1$. Nous allons minorer 
$\textrm{dim}_H \pi_D\left({\mathcal G} \cap Z\right)$. Commen\c cons par observer que tous les graphes de $\mathcal G$ intersectent $Z$.

\begin{lem}\label{Lem} 
L'intersection $\Gamma_\gamma \cap Z$ est non vide
pour tout  $\gamma\in {\mathcal O}\left(D,B\right)$.  En particulier, $\Gamma_\omega \cap Z$ est non vide et discret pour tout $\omega\in \Sigma_m^+$.
\end{lem}

\proof 
L'espace ${\mathcal O}\left(D,B\right)$ est convexe, on le munit de la topologie de la convergence uniforme locale. Pour tout $\gamma\in {\mathcal O}\left(D,B\right)$,
le sous-ensemble analytique $\pi_D\left(\Gamma_\gamma\cap Z\right)$ de $D$ est de dimension nulle car il est
relativement compact dans $D$. Ainsi $\Gamma_\gamma \cap Z$ est discret (ou vide). L'ensemble $\{\gamma\in {\mathcal O}\left(D,B\right)\; :\; \Gamma_\gamma \cap Z \ne \emptyset\}$ est clairement ferm\'e dans  ${\mathcal O}\left(D,B\right)$, le lemme d'Hurwitz montre qu'il est ouvert et la conclusion
s'ensuit. \fin

\begin{prop}\label{PropZ} Pour toute hypersurface irreductible $Z \subset D\times B$, non
 r\'eduite \`a une fibre de $\pi_D$ et telle que $\pi_D(Z) \Subset D$  on a l'estimation suivante :
$\textrm{dim}_H \pi_D\left({\mathcal G} \cap Z\right) \ge \frac{a}{A} \left(\frac{\ln m}{A}\right) - (2k-2)$.
\end{prop}

\proof On suppose que $\frac{\ln m}{A} > \frac{A}{a} (2k-2)$ car sinon il n'y a rien \`a d\'emontrer. Voyons d'abord comment se ramener au cas o\`u $Z$ est lisse et transverse aux fibres de $\pi_D$. Rappelons que $Z$ est un sous-ensemble analytique de dimension complexe $k$ dans $D\times B$. Soit $S$ l'ensemble constitu\'e des points
singuliers de $Z$ ainsi que de ses points r\'eguliers $p$ pour lesquels l'espace tangent $T_p(Z)$ est confondu avec la fibre
$\pi_D^{-1}\left(\pi_D(p)\right)$. Nous allons montrer qu'il existe $\tilde{\omega}\in \Sigma_m^+$ tel que 
$\Gamma_{\tilde{\omega}} \cap Z$ n'est pas contenu dans   $S$.\\
Comme $S$ est un sous-ensemble analytique strict de $Z$ (par hypoth\`ese $Z$ n'est pas une fibre de $\pi_D$) on a $\textrm{dim}_H \left(S\right) \le 2k-2$. Supposons que $\Gamma_{{\omega}} \cap Z \subset  S$ pour tout 
$\omega \in \Sigma_m^+$ et fixons $\la_0 \in D$. Comme $\Gamma_\omega \cap Z \ne \emptyset$ on a
$H_{\la_0}\left(\Gamma_\omega\right) =H_{\la_0}\left(\Gamma_\omega \cap Z\right)$ pour tout $\omega\in \Sigma_m^+$ et donc $\left({\mathcal G}\right)_{\la_0}= H_{\la_0}\left({\mathcal G}\right) = \cup_{\omega \in \Sigma_m^+} H_{\la_0}
\left(\Gamma_\omega\right) \subset H_{\la_0} \left(S\right)$. D'apr\`es la Proposition \ref{PropCantGr}, cela entra\^{i}ne que $\frac{\ln m}{A} \le \textrm{dim}_H \left(\left({\mathcal G}\right)_{\la_0}\right) \le \frac{A}{a} \textrm{dim}_H \left(S\right) \le \frac{A}{a} (2k-2)$ ce qui est exclu. Il existe donc $\tilde{\omega}\in \Sigma_m^+$ et $\la_1 \in D_{r_0}$
tels que $(\la_1,\tilde{\omega}(\la_1)) \in Z\setminus S$.
Etant donn\'e  un voisinage $V$ de $(\la_1,\tilde{\omega}(\la_1))$ dans $Z\setminus S$ on voit, gr\^{a}ce au lemme d'Hurwitz, que si $p$ est assez grand alors $\Gamma_\omega \cap V \ne \emptyset$ pour tout $\omega \in \Sigma_m^+$ v\'erifiant $\omega_j={\tilde \omega}_j$ pour $j\le p$.\\

Nous obtiendrons l'estimation annonc\'ee en rempla\c cant $\mathcal G$ par  $\widetilde{\mathcal G}_p:=\cup_{\omega\in \widetilde{C}_p}\Gamma_\omega$ o\`u
$\widetilde{C}_p:=\{\omega \in \Sigma_m^+\; :\; \omega_j={\tilde \omega}_j \;\textrm{si}\;j\le p \}$ . Observons que $(\widetilde{\mathcal G}_p)_{\la} =G_{\tilde{\omega}_0\cdots \tilde{\omega}_p} \left(\left({\mathcal G}\right)_{\la}\right)$ pour tout $\la \in D$.
Comme $G_{\tilde{\omega}_0\cdots \tilde{\omega}_p} $ est bi-lipschitzienne sur $\{\la_0\}\times B$, il  r\'esulte de la Proposition \ref{PropCantGr} que $\textrm{dim}_H ((\widetilde{\mathcal G}_p)_{\la_0}) \ge \frac{\ln m}{A}$ puis, comme
$(\widetilde{\mathcal G}_p)_{\la_0}=H_{\la_0}(\widetilde{\mathcal G}_p) \subset H_{\la_0}(\widetilde{\mathcal G}_p \cap Z)$, que
$\frac{\ln m}{A}\le \frac{A}{a} \textrm{dim}_H (\widetilde{\mathcal G}_p \cap Z)$. Posons
$E_p:=\widetilde{\mathcal G}_p \cap Z$; il nous reste pour conclure \`a justifier la majoration
$\textrm{dim}_H \left(E_p\right) \le \textrm{dim}_H \left(\pi_D\left({\mathcal G}\cap Z\right)\right) + (2k-2)$.\\
On peut supposer que $V=\{(l(z),z)\; :\; z\in U\}$ o\`u $U$ est un voisinage assez petit de $\tilde{\omega}(\la_1)$ et $l \in  {\mathcal O}(U,\C)$
v\'erifie $\frac{\partial l}{\partial z_1} \ne 0$ sur $U$. Alors, quitte \`a diminuer $U$, l'application $\psi(\la,z):=(\la,l(z),z_2,\cdots,z_k)$ induit un 
biholomorphisme sur un voisinage $\Omega_{\la_1}\times U$ de $(\la_1,\tilde{\omega}(\la_1))$ et, puisque 
$E_p\Subset \Omega_{\la_1}\times U$ pour $p$ est assez grand, il vient $\textrm{dim}_H \left(E_p\right) =\textrm{dim}_H \left(\psi(E_p)\right)$.
On termine en remarquant que $\psi(E_p)\subset \{(\la,\la)\; :\; \la \in \pi_D(E_p)\} \times \C^{k-1}$ et $\textrm{dim}_H \left(\{(\la,\la)\; :\; \la \in \pi_D(E_p)\} \times \C^{k-1}\right)
\le \textrm{dim}_H \left(\{(\la,\la)\; :\; \la \in \pi_D(E_p)\}\right) +    (2k-2)
\le \textrm{dim}_H \left(\pi_D\left({\mathcal G}\cap Z\right)\right) + (2k-2)$. \fin

\section{Contractions issues d'une perturbation dans ${\mathcal H}_d\left(\Pj^k\right)$}

Consid\'erons  $f\in {\mathcal H}_d\left(\Pj^k\right)$ dont les exposants de Lyapounov $\chi_1\le\cdots \le\chi_k$ relatifs \`a sa mesure d'entropie maximale $\mu_f$ ne satisfont aucune relation de r\'esonnance (i.e. $\alpha_1\chi_1+\cdots+\alpha_k\chi_k\ne \chi_j$ pour tout $1\le j\le k$ et tout $\alpha \in \N^k$ tel que 
$\alpha_1+\cdots+\alpha_k \ge 2$). Notre objectif est d'associer une lamination du type de celles \'etudi\'ees \`a la section pr\'ec\'edente \`a toute famille holomorphe $F:D\times \Pj^k \to D\times \Pj^k$ telle que $F(0,\cdot)=f$. Nous construirons pour cela des branches inverses it\'er\'ees de $f$ dont les distorsions sont contr\^{o}l\'ees puis les prolongerons en des branches inverses de $F$. Ceci repose essentiellement sur l'application \`a $f$ d'une m\'ethode de lin\'earisation le long des orbites \'etablie dans \cite{BDM} et dont nous allons commencer par rappeler le principe.\\

 Soit $\textrm{O}:=\{\hat{x}:=(x_n)_{n\in \Z}\; :\; f(x_n)=x_{n+1}\}$ l'espace des orbites. On note $\pi$ la projection ${\hat x} \mapsto x_0$ puis 
 $\hat f$ le d\'ecalage \`a gauche ($\pi \circ {\hat f} =f\circ \pi$) et $\tau$ son inverse. Pour tout $E\subset \Pj^k$ on posera
 ${\widehat E}:=\pi^{-1}\left( E\right)$. Il existe une unique mesure de probabilit\'e $\nu$ sur 
 $\textrm{O}$ telle que $\pi_{\star} \nu =\mu_f$, cette mesure est m\'elangeante. On travaillera dans l'espace $X:=\{{\hat x}\in \textrm{O}\; :\;
x_n \notin {\mathcal C}_f, \forall n\in \Z\}$ o\`u $C_f$ d\'esigne l'ensemble critique de $f$. Comme $\mu_f\left({\mathcal C}_f\right) =0$, l'ensemble $X$ est de mesure pleine pour $\nu$. Pour tout ${\hat x}\in X$, on note $f_{\hat x}^{-n}$ la branche inverse de $f^n$ d\'efinie au voisinage de $x_0$ et envoyant $x_0$ sur $x_{-n}$. Rappelons qu'une fonction $\alpha : \textrm{O}
\to ]0,1]$ est dite $\epsilon$-lente si $\alpha \circ \tau \ge e^{-\epsilon} \alpha$. \\
On note $d(\cdot\;,\cdot)$ la distance induite sur $\Pj^k$ par la m\'etrique de Fubini-Study, $ B_x(r)\subset \Pj^k$ la boule centr\'ee en $x$ et de rayon $r$ pour cette distance et $B(r)$ la boule $\{\Vert z\Vert <r\} \subset \C^k$ pour la distance hermitienne usuelle.\\
 
 Le th\'eor\`eme 1.4 de \cite{BDM} stipule que, pour $0<\epsilon \ll \chi_1$, il existe des fonctions $\epsilon$-lentes $r_\epsilon$, $t_{\epsilon}$, $1/\beta_{\epsilon}
 :X\to ]0,1]$, une constante $0<\alpha\le 1$, des applications injectives $S_{\hat x}$
 et des applications \emph{lin\'eaires} $R_{\hat x}^n$ telles que le diagramme 
\begin{center}
\begin{tikzpicture}
\matrix (m) [matrix of math nodes, row sep=3.5em,
column sep=3.5em, text height=1.5ex, text depth=0.25ex]
{ B_{x_0} (r_\epsilon (\hat x)) & f_{\hat x}^{-n} ( B_{x_0} (r_\epsilon (\hat x))) \\
  B(t_\epsilon (\hat x)) & B (t_\epsilon (\tau^n (\hat x)))\\};
	\path[->,font=\scriptsize]
		(m-1-1) edge node[left] {$ S_{\hat x} $} (m-2-1)
		(m-1-2) edge node[right] {$ S_{\tau^n (\hat x)} $} (m-2-2)
		(m-1-1) edge node[above] {$ f_{\hat x}^{-n} $} (m-1-2)
		(m-2-1) edge node[above] {$ R^n_{\hat x}$} (m-2-2);
\end{tikzpicture}
\end{center}
 commute pour
 tout $n\in\N$ et tout $\hat x \in X$.
 
 Soulignons que $S_{\hat x}(x_0)=0$ et que la lin\'earit\'e des applications $R_{\hat x}^n$ d\'ecoule de l'absence de r\'esonnance sur les exposants de Lyapounov de $f$. En outre, les applications $S_{\hat x}$ et $R_{\hat x}^n$ satisfont les estimations suivantes
 \begin{eqnarray}\label{aa}
 e^{-n\chi_k} \Vert z\Vert \le \Vert R_{\hat x}^n (z)\Vert \le  e^{-n\chi_1} \Vert z\Vert
 \end{eqnarray}
  \begin{eqnarray}\label{bb}
 \alpha d(p,q) \le \Vert S_{\tau^n(\hat x)}(p) -  S_{\tau^n(\hat x)}(q) \Vert \le \beta_{\epsilon}\left(\tau^n({\hat x})\right) d(p,q)
 \end{eqnarray}
   \begin{eqnarray}\label{cc}
e^{-n (\chi_k+\epsilon)} \frac{\alpha}{\beta_\epsilon({\hat x})} d(p,q) \le d(f_{{\hat x}}^{-n} (p),f_{{\hat x}}^{-n} (q)) \le  
e^{-n \chi_1} \frac{\beta_\epsilon({\hat x})}{\alpha} d(p,q)
 \end{eqnarray}
l'assertion (\ref{cc}) \'etant une cons\'equence directe de la commutativit\'e du diagramme et des assertions (\ref{aa}) et (\ref{bb}).\\

Posons $\rho_\epsilon:=\alpha\frac{r_\epsilon}{\beta_\epsilon}$ (c'est une fonction $2\epsilon$-lente telle que 
$0<\rho_\epsilon \le r_\epsilon$) puis, pour tout $0<t\le 1$ d\'efinissons les ensembles $E_{\hat x}^{-n}(t) \subset \widetilde{E}_{\hat x}^{-n}(t)$ par  
$E_{\hat x}^{-n}(t):=f_{\hat x}^{-n}\left(
B_{x_0}\left(t\rho_\epsilon({\hat x})\right)\right)$
et $\widetilde{E}_{\hat x}^{-n}(t):=f_{\hat x}^{-n}\left(
B_{x_0}\left(t r_\epsilon({\hat x})\right)\right)$.
L'objet du lemme suivant est de  pr\'eciser la g\'eom\'etrie des ensembles $E_{\hat x}^{-n}(t)$.
 \begin{lem}\label{lemConv}
Pour  tout $0<t\le 1$ et toute paire de points $ p,q$ de $E_{\hat x}^{-n}(t)$, il existe un chemin lisse joignant $p$ \`a $q$  dans $\widetilde{E}_{\hat x}^{-n}(t)$
  et dont la longueur est major\'ee par $ \frac{\beta_\epsilon({\hat x})}{\alpha} e^{n\epsilon} d(p,q)$.
  \end{lem}
  
 \proof Si $p,q\in E_{\hat x}^{-n}(t)$ alors $f^n(p),f^n(q)\in B_{x_0}\left(t\rho_\epsilon(\hat x)\right)$ et 
 $R_{\hat x}^n\circ S_{\hat x} \left(f^n(p)\right) = S_{\tau^n ({\hat x})}(p)$, $R_{\hat x}^n\circ S_{\hat x} \left(f^n(q)\right) = S_{\tau^n ({\hat x})}(q)$. On voit ainsi, en utilisant (\ref{bb}), que
 \begin{eqnarray*}
S_{\tau^n ({\hat x})}(p), S_{\tau^n ({\hat x})}(q) \in R_{\hat x}^n\left[ S_{\hat x} \left(B_{x_0}\left(t \rho_\epsilon(\hat x)\right)\right)\right] \subset R_{\hat x}^n\left[B( t \alpha   r_\epsilon ({\hat x}))\right].
 \end{eqnarray*}
 V\'erifions \'egalement que 
 \begin{eqnarray*}
 R_{\hat x}^n\left[B(t\alpha r_\epsilon ({\hat x}))\right] \subset  S_{\tau^n ({\hat x})}\left[ B_{x_{-n}} \left(t r_\epsilon (\tau^n({\hat x}))\right)  \right].
 \end{eqnarray*}
Cela r\'esulte des inclusions   
$ R_{\hat x}^n\left[B(t \alpha r_\epsilon ({\hat x}))\right] \subset B\left(t \alpha r_\epsilon ({\hat x}) e^{-n\chi_1}\right) \subset B(t \alpha  r_\epsilon ({\hat x})e^{-n\epsilon})$ et 
$B(t \alpha  r_\epsilon ({\hat x})e^{-n\epsilon}) \subset B(t \alpha r_\epsilon ({\tau^n(\hat x})))
\subset  S_{\tau^n ({\hat x})}\left[ B_{x_{-n}} \left(t r_\epsilon (\tau^n({\hat x}))\right)  \right]$
qui se d\'eduisent respectivement de 
 (\ref{aa}) et  (\ref{bb}).
 
 On peut donc consid\'erer l'image par $\left(S_{\tau^n ({\hat x})}\right)^{-1}$ du segment joignant  $S_{\tau^n ({\hat x})}(p)$ \`a $ S_{\tau^n ({\hat x})}(q)$ dans le \emph{convexe} $R_{\hat x}^n\left[B(t\alpha r_\epsilon ({\hat x}))\right]$. C'est un chemin joignant $p$ \`a $q$. En utilisant
(\ref{bb}), on voit que sa longueur est major\'ee par $\frac{1}{\alpha} \beta_\epsilon (\tau^n({\hat x})) d(p,q) \le 
\frac{\beta_\epsilon ({\hat x})}{\alpha} e^{n\epsilon} d(p,q)$. Enfin, ce chemin est contenu dans $\widetilde{E}_{\hat x}^{-n}(t)$ car 
$\left(S_{\tau^n ({\hat x})}\right)^{-1}\left[R_{\hat x}^n\left[B(t\alpha r_\epsilon ({\hat x}))\right]\right] \subset
\left(S_{\tau^n ({\hat x})}\right)^{-1} \circ R_{\hat x}^n \circ S_{\hat x}\left[B_{x_0}(t r_\epsilon ({\hat x}))\right] = \widetilde{E}_{\hat x}^{-n}(t)$. \fin

La Proposition suivante est obtenue en combinant les r\'esultats d\'ecrits ci-dessus avec un argument classique d\^{u} \`a Briend-Duval.

\begin{prop}\label{PropBDBDM} Soit $f\in {\mathcal H}_{d\ge 2} \left(\Pj^k \right)$. On suppose que   les exposants de Lyapounov  $\chi_1\le \cdots \le \chi_k$ de $f$ par rapport \`a sa mesure d'entropie maximale $\mu_f$ ne satisfont aucune relation de r\'esonnance. Alors, pour tout ouvert $\Omega \subset \Pj^k$ tel que $\mu_f\left( \Omega\right)>0$ et tout $0<\epsilon \ll 1$, il existe des constantes $r_0,C,K>0$ (d\'ependant de $\epsilon$) ainsi qu'une boule $A\subset \Omega$ charg\'ee par $\mu_f$ telles que, pour tout $n$ assez grand, $f^n$ admet $m\ge Cd^{kn}$
branches inverses $g_1,\cdots, g_m$ d\'efinies sur des boules $B_{x_0^j} (Kr_0)$  v\'erifiant
$B_{x_0^j} (Kr_0) \supset B_{x_0^j} (r_0) \supset A$ et
\begin{itemize}
\item[1)] $g_j(A) \Subset A$,
\item[2)] $e^{-n(\chi_k +\epsilon)} \frac{1}{K} d(p,q) \le d\left(g_j(p),g_j(q)\right) \le e^{-n\chi_1} K d(p,q)$ sur $B_{x_0^j}(Kr_0)$,
\item[3)] deux points quelconques $ p,q$ de $g_j(B_{x_0^j}(r_0))$ sont reli\'es par un chemin lisse dans $g_j(B_{x_0^j}(Kr_0))$
dont la longueur est major\'ee par $ K e^{n\epsilon} d(p,q)$.
\end{itemize}
\end{prop}

\proof 
Rappelons que $\rho_\epsilon =\alpha \frac{r_\epsilon}{\beta_\epsilon}$ et d\'efinissons une partie $\widehat H$ de $X$ par
$\widehat{H}:=\{\hat x \in X\; :\; r_\epsilon(\hat x) \ge K r_0, \frac{\alpha}{\beta_\epsilon(\hat x)} \ge \frac{1}{K}\}$. Choisissons
 $K\gg 1$ puis $0<r_0\ll 1$ pour que $\nu\left( {\widehat H} \cap {\widehat \Omega}\right) >0$.\\
Par un argument de recouvrement, on trouve une boule $A_r$ de rayon $r>0$ dans $\Omega$ et $0<\gamma \ll r$ tels que $r+\gamma
<\frac{r_0}{2}$ et $\nu\left( {\widehat H} \cap {\widehat A_{r+\gamma}}\right) >0$. Observons que si $\hat x \in {\widehat H} \cap {\widehat A_{r+\gamma}}$ alors $A_{r+\gamma} \subset B_{x_0}(2(r+\gamma)) \Subset B_{x_0} (r_0)\subset B_{x_0} (Kr_0)$ et qu'alors $f_{\hat x}^{-n}$ est bien d\'efinie sur $A_{r+\gamma}$. On d\'efinit alors $C_n$
par $C_n:=\{\hat x \in {\widehat H} \cap {\widehat A_{r+\gamma}}\; : \; f_{\hat x}^{-n}\left( A_{r+\gamma}\right) \cap A_r \ne \emptyset\}$.\\
\`A chaque $\hat x\in C_n$ correspond une branche inverse $g:=f_{\hat x}^{-n}$ de $f^n$ d\' efinie sur une boule $B_{x_0}(Kr_0) \supset A_{r+\gamma}$. Il r\'esulte du th\'eor\`eme 1.4 de \cite{BDM} (dont le contenu est rappel\'e ci-dessus) et de la d\'efinition de $\widehat H$ que $g$ satisfait les assertions (2) et (3). Plus pr\'ecis\'ement, l'assertion $(2)$ est une reformulation de  (\ref{cc}) tenant compte de $\frac{\alpha}{\beta_\epsilon(\hat x)} \ge \frac{1}{K}$ et l'assertion (3) s'obtient en appliquant le  Lemme \ref{lemConv} avec $t:=\frac{r_0}{\rho_\epsilon (\hat x)}$. Il est utile  d'observer que $t=\frac{r_0}{r_\epsilon (\hat x)} \frac{r_\epsilon (\hat x)}{\rho_\epsilon (\hat x)}
\le \frac{1}{K} \frac{\beta_\epsilon (\hat x)}{\alpha}\le 1$
et que $t{r_\epsilon (\hat x)} =r_0\frac{r_\epsilon (\hat x)}{\rho_\epsilon (\hat x)}= r_0
\frac{\beta_\epsilon (\hat x)}{\alpha}\le Kr_0$.

En particulier, puisque $g\left(A_{r+\gamma}\right) \cap A_r \ne \emptyset$ et $\textrm{diam}\;g\left(A_{r+\gamma}\right) \le 2(r+\gamma) e^{-n\chi_1} K$, on voit que $g\left(A_{r+\gamma}\right) \Subset A_{r+\gamma}$ pourvu que $e^{n\chi_1}>2\frac{r+\gamma}{\gamma} K $. On prendra donc $A:=A_{r+\gamma}$.\\
 Il reste \`a \'etablir l'assertion $(1)$ et, pour cela, \`a minorer le nombre de branches inverses distinctes au-dessus de $A$ associ\'ees aux  $\hat x\in C_n$. Notons $m$ ce nombre.
Observons  que deux \'el\'ements distincts de $C_n$ 
 donnent lieu \`a  des images $f_{\hat x}^{-n} (A)$ disjointes ou confondues.
Nous  allons combiner le caract\`ere m\'elangeant de $\nu$ et
le fait que le jacobien de $\mu$ est  constant \'egal \`a $d^k$.
La premiere propri\'et\'e implique que
\[
\nu (\tau^n (\widehat H \cap \widehat {A}_{}) \cap \widehat {A}_{r})
\to
\nu (\widehat H \cap \widehat {A}_{}) \cdot \nu (\widehat {A}_{r})
\]
pour $n\to\infty$.
En utilisant le fait que $f^* \mu = d^k \mu$ on a donc
\[
\begin{aligned}
m\; \mu(A_{}) d^{-kn}
& \geq
\mu
(\cup_{C_n} f_{\hat x}^{-n} (A_{}))\\
& \geq
\nu
(\tau^n (\widehat H \cap \widehat {A}_{}) \cap \widehat {A}_{r})
\geq
\nu (\widehat H \cap \widehat {A}_{}) \nu (\widehat {A}_r)/2 >0
\end{aligned}
\]
pour $n$ assez grand, ce qui donne l'estimation annonc\'ee.\fin\\

Consid\'erons maintenant une famille holomorphe $F:D\times \Pj^k \to D\times \Pj^k$. Notons $f_\la:=F(\la,\cdot)$ l'endomorphisme de $\Pj^k$ correspondant au param\`etre $\la$ et $\mu_\la$ sa mesure d'entropie maximale. Soit $\Omega \subset \Pj^k$ un ouvert charg\'e par $\mu_0$ et $g_{0,1},\cdots,g_{0,m}$ les branches inverses de $f_0^n$ fournies par la Proposition \ref{PropBDBDM}. Nous montrerons que, quitte \` a diminuer l\'eg\`erement $A$, les $g_{0,j}$
se prolongent en des biholomorphismes $G_j : D_{r(n)} \times {A} \to G_j(D_{r(n)} \times { A})$ qui h\'eritent des $g_{0,j}$ leurs propri\'et\'es contractantes. C'est ici que l'assertion $(3)$ de la Proposition \ref{PropBDBDM} et  implicitement le th\'eor\`eme
de lin\'earisation de \cite{BDM} jouent un r\^{o}le crucial. Cela conduira \`a la proposition suivante.

\begin{prop}\label{PropJC} Soit $F:D\times \Pj^k \to D\times \Pj^k$ une famille holomorphe d'endomorphismes de degr\'e $d\ge 2$. 
On suppose que   les exposants de Lyapounov  $\chi_1\le \cdots \le\chi_k$ de l'endomorphisme
$f_0:=F(0,\cdot)$ par rapport \`a sa mesure d'entropie maximale $\mu_0$ ne satisfont aucune relation de r\'esonnance. Alors, pour tout ouvert $\Omega \subset \Pj^k$ tel que $\mu_0\left( \Omega\right)>0$ et tout $0<\epsilon \ll 1$, il existe une constante $C>0$ (d\'ependant de $\epsilon$) ainsi qu'une boule ${ A}\subset \Omega$ charg\'ee par $\mu_0$ telles que, pour tout $n$ assez grand, $F^n$ admet $m\ge Cd^{kn}$
branches inverses $G_1,\cdots, G_m$ d\'efinies au voisinage de  $\overline{D_{r(n)}}\times \overline{ A}$  et v\'erifiant
\begin{itemize}
\item[1)] $G_j(D_{r(n)} \times { A}) \subset D_{r(n)} \times { A}$,\;$\forall 1\le j\le m$,
\item[2)] $\textrm{dist} \left(G_j(D_{r(n)} \times { A}), G_k (D_{r(n)} \times { A})\right) >0,\; \forall 1\le j\ne k\le m$,
\item[3)] $e^{-n(\chi_k +3\epsilon)}  d(z,z') \le d\left(G_j(\la,z),G_j(\la,z')\right) \le e^{-n(\chi_1-\epsilon)}  d(z,z')$ sur $D_{r(n)} \times { A}$.
\end{itemize}
\end{prop}

\proof Rappelons que les branches $g_{0,j}$ donn\'ees par la Proposition \ref{PropBDBDM} sont toutes d\'efinies sur une m\^{e}me boule $\widetilde A$ (contenue dans $\Omega$) et qu'elles induisent des biholomorphismes $g_{0,j} : B_{x_0^j}(r_0) \to g_{0,j}(B_{x_0^j}(r_0))$ o\`u ${\widetilde A}\Subset B(x_0^j,r_0)$. 

Soit $A$ une boule relativement compacte dans $\widetilde A$, obtenue en diminuant l\'eg\`erement le rayon de $\widetilde A$ de fa\c con \`a  ce que $ A$ reste charg\'ee par $\mu_0$. Posons $A_j:=
g_{0,j}( \widetilde A)$, c'est
un ouvert \`a bord r\'egulier contenu dans $g_{0,j}(B_{x_0^j}(r_0))$ et $f_0^n(bA_j) = b\left(f_0^n(A_j)\right)$. Les $f_\la^n$ \'etant continues et ouvertes, on a 
$b\left(f_\la^n(A_j)\right)\subset f_\la^n (bA_j)$
et l'on d\'eduit alors de ${ A} \Subset f_0^n\left(A_j\right)$ 
que ${ A} \Subset f_\la^n\left(A_j\right)$ pour $\la \in D_{r(n)}$ pourvu que $r(n)$ soit assez petit. Cela se traduit par 
$$D_{r(n)}\times { A}\subset
F^n\left( D_{r(n)}\times { A_j}\right).$$
 Comme $f_0^n$ ne branche pas au voisinage de $\overline {A_j}$, on peut diminuer $r(n)$ de fa\c con \`a ce que $F^n$ ne branche pas sur 
$D_{r(n)}\times { A_j}$. Soit $U_j$ une composante connexe de $(F^n\vert_{D_{r(n)}\times{ A_j}})^{-1} ( D_{r(n)}\times{ A})$, alors l'application 
$F^n : U_j\to D_{r(n)}\times{ A}$ est holomorphe propre et ne branche pas. Puisque $D_{r(n)}\times{ A}$ est simplement connexe c'est un biholomorphisme,
on d\'efinit $G_j$ comme \'etant son inverse. Quitte \`a diminuer l\'eg\`erement $r(n)$ et le rayon de $ A$, on pourra supposer que $G_j$ est d\'efinie au voisinage de l'adh\'erence de ${D_{r(n)}\times{ A}}$.\\

Voyons maintenant que les $G_j$ satisfont les estimations annonc\'ees. Il faudra pour cela r\'eduire \`a plusieurs reprises $r(n)$, ce que nous ferons sans le pr\'eciser. 
Les $G_j$ sont de la forme $G_j(\la,z)=\left(\la,g_{\la,j} (z)\right)$ o\`u $f_\la^n\circ g_{\la,j}=id$.
 Puisque $g_{0,j}({ A}) \Subset { A}$ et 
$d(g_{0,j}({ A}),g_{0,k}({ A}))>0$ si $j\ne k$ on a $g_{\la,j}({ A}) \Subset { A}$ et 
$\textrm{dist}(G_j(D_{r(n)}\times{ A}),G_k(D_{r(n)}\times{ A}))>0$ si $j\ne k$.\\
\'Etablissons l'estimation interm\'ediaire 
\begin{eqnarray}\label{int}
e^{-n(\chi_k +2\epsilon)} \frac{1}{2K^2} d(z,z') \le d\left(G_j(\la,z),G_j(\la,z')\right) \le e^{-n\chi_1} 2K d(z,z')\;\textrm{sur}\; D_{r(n)} \times { A}.
\end{eqnarray}
Par la Proposition \ref{PropBDBDM}, on a $d\left(g_{0,j}(z),g_{0,k}(z')\right)\le  K e^{-n\chi_1}  d(z,z') $ pour $z,z'\in B_{x_0^j}(r_0)$. Il s'ensuit que $\Vert d_z G_j(0,z)\Vert \le e^{-n\chi_1} K$
sur $ A$ puis que $\Vert d_z G_j(\la,z)\Vert \le 2e^{-n\chi_1} K$ sur $D_{r(n)}\times { A}$. Ceci entra\^{i}ne la majoration dans (\ref{int}).

Par la Proposition \ref{PropBDBDM}, on a $d\left(f_0^n(z),f_0^n(z')\right)\le K e^{n(\chi_k+\epsilon)} d(z,z')$ pour tous $z,z'\in g_{0,j}(B_{x_0^j}(r_0))$. Il s'ensuit que  $\Vert d_z f_0^n(z)\Vert\le K e^{n(\chi_k+\epsilon)}$
sur $g_{0,j}(B_{x_0^j}(r_0))$ et donc  que $\Vert d_z F^n(\la,z)\Vert\le 2K e^{n(\chi_k+\epsilon)}$ sur $D_{r(n)}\times g_{0,j}(B_{x_0^j}(r_0))$.
Puisque $g_{0,j}({\widetilde A}) \Subset g_{0,j}(B_{x_0^j}(r_0))$ on a aussi $g_{\la,j}({ A}) \Subset g_{0,j}(B_{x_0^j}(r_0))$ pour $\la \in D_{r(n)}$. Ainsi, compte tenu de la troisi\`eme assertion de la Proposition \ref{PropBDBDM}, 
si $(\la,z)$ et $(\la,z')$ sont dans $D_{r(n)}\times { A}$ alors les points $g_{\la,j}(z)$ et $g_{\la,j}(z')$ sont reli\'es par un chemin dont la longueur est major\'ee par $Ke^{n\epsilon} 
d(g_{\la,j}(z) , g_{\la,j}(z')) $ et sur lequel $\Vert d_z F^n(\la,z)\Vert\le 2K e^{n(\chi_k+\epsilon)}$.  La minoration dans  (\ref{int}) s'ensuit car $d(z,z') =
d(F^n(\la,g_{\la,j}(z)),F^n(\la,g_{\la,j}(z'))) \le  2K e^{n(\chi_k+\epsilon)} Ke^{n\epsilon} 
d(g_{\la,j}(z) , g_{\la,j}(z')) $. L'assertion (3) d\'ecoule de (\ref{int}) lorsque $n$ est assez grand.\fin

\section{Dimension de Hausdorff du lieu de bifurcation}

Consid\'erons  une   famille holomorphe $F:M\times \Pj^k \to M\times \Pj^k$ d'endomorphismes de $\Pj^k$, notre objectif
est d'obtenir des minorations locales de la dimension de Hausdorff de son lieu de bifurcation $\Bif$. Nous allons pour cela utiliser les r\'esultats des sections pr\'ec\'edentes ainsi que la caract\'erisation suivante du lieu de bifurcation (voir le Th\'eor\`eme 1.6
de \cite{BBD}).

\begin{thm}\label{thMis} L'ensemble des  param\` etres du type Misiurewicz est une partie dense du lieu de bifurcation.
\end{thm}

\begin{defn}\label{DefiMisiu}
On dit que  $\lambda_0 \in D$ est un \emph{param\`etre du type Misiurewicz } si il existe une application holomorphe $\gamma $ d\'efinie sur un voisinage de  $\lambda_0$ et \`a  valeurs dans $\Pj^{k}$ telle que :
\begin{itemize}
\item[1)]  $\gamma(\lambda_0)\in J_{\lambda_0}$ et $\gamma(\lambda)$ est $p_0$-periodique r\'epulsif pour $f_{\lambda}$ et un certain $p_0 \geq 1$,
\item[2)] $(\lambda_0,\gamma(\lambda_0))\in F^{n_0}(C_F)$ pour un certain $n_0 \geq 1$,
\item[3)] le graphe $\Gamma_\gamma$ de $\gamma$ n'est pas contenu dans  $F^{n_0}(C_F)$.
\end{itemize}
\end{defn}

Rappelons que  l'ensemble de Julia $J_\la$
de $f_\la:=F(\la,\cdot)$ est, par d\'efinition, le support de  la mesure d'entropie maximale $\mu_\la$ de $f_\la$ et que $C_F$ d\'esigne l'ensemble critique
de l'application $F:M\times \Pj^k \to M\times \Pj^k$.\\

Le principal r\'esultat de cette section est le suivant.

\begin{prop}\label{PropMis} Soit $F:D\times \Pj^k \to D\times \Pj^k$ une famille holomorphe d'endomorphismes de degr\'e $d\ge 2$. 
On suppose que $0$ est un param\`etre Misiurewicz et que  les exposants de Lyapounov  $\chi_1\le \cdots \le \chi_k$ de $f_0:=F(0,\cdot)$ par rapport \`a sa mesure d'entropie maximale $\mu_0$ ne satisfont aucune relation de r\'esonnance. On a alors l'estimation suivante : $\liminf_{r\to 0} \textrm{dim}_H\left(D_r\cap {\Bif} \right)\ge \frac{\chi_1}{\chi_k} (\frac{k\ln d}{\chi_k})-(2k-2)$.
\end{prop}

Comme cons\'equence nous obtenons que la dimension de Hausdorff du lieu de bifurcation au voisinage d'un exemple de Latt\`es est  maximale; c'est le th\'eor\`eme \ref{thP}.\\

\noindent\textsc{D\'emonstration du th\'eor\`eme  \ref{thP}:} Notons $\chi_1(\la)\le \cdots\le \chi_k(\la)$ les exposants de Lyapounov de $f_\la$ et rappelons que $\chi_1(\la) \ge \frac{\ln d}{2}$. Que $0$ soit dans le lieu de bifurcation r\'esulte directement de la caract\'erisation des exemples de Latt\`es par la minimalit\'e de leurs exposants  et de celle de ${\Bif}$ comme \'etant le support de 
$dd^c (\chi_1(\la)+\cdots+\chi_k(\la))$ (voir  les Th\'eor\`emes \ref{THBDL} et \ref{THBBD}), pour plus de d\'etails nous renvoyons au Th\'eor\`eme 6.3 de \cite{BBD}. En vertu du Th\'eor\`eme \ref{thMis},
on peut approcher $0$ par des param\`etres Misiurewicz $\la_n$. La somme des  exposants de Lyapounov d\'ependant contin\^{u}ment du param\` etre (voir \cite[Theorem 2.47]{DS} )  il vient $\lim_n\chi_k(\la_n) =\chi_k(0)=\frac{\ln d}{2}$ et  $\lim_n\chi_1(\la_n)=\frac{\ln d}{2}$. En particulier les $\chi_j(\la_n)$ ne satisfont aucune relation de r\'esonnance pour $n$ assez grand. La conclusion d\'ecoule alors de la Proposition \ref{PropMis}.\fin

\begin{rem}
En dimension $k=1$, l'exposant de Lyapounov $\chi(\la)$ est reli\'e \`a la dimension de Hausdorff ${\textrm dim}_H \mu_\la$ de la mesure $\mu_\la$ par la formule $\ln d= \chi(\la) {\textrm dim}_H \mu_\la$. La Proposition \ref{PropMis} montre alors que pour toute famille holomorphe $F:D\times \Pj^1\to D\times \Pj^1$ de fractions rationnelles de degr\'e $d\ge 2$ telle que $0$ appartienne au lieu de bifurcation, on a
$\liminf_{r\to 0} \textrm{dim}_H\left(D_r\cap {\Bif} \right)\ge  {\textrm dim}_H \mu_\la$. Notons cependant que cette minoration est loin d'\^{e}tre optimale car, comme l'a montr\'e M. Shishikura, la dimension de Hausdorff du bord de l'ensemble de Mandelbrot est \'egale \`a deux \cite{Shi2}.
\end{rem}

\noindent\textsc{D\'emonstration de la proposition \ref{PropMis}:} Observons que l'on peut remplacer la famille $F$ par une it\'er\'ee $F^q$ car le lieu de bifurcation ainsi que  la quantit\'e   $\frac{\chi_1}{\chi_k} (\frac{k\ln d}{\chi_k})$ apparaissant dans la minoration de $\textrm{dim}_H\left(D_r\cap {\Bif}\right)$  restent alors  inchang\'es. La d\'emonstration reprend le principe de celle mise au point par le second auteur pour g\'en\'eraliser le Th\'eor\`eme
\ref{thMis} aux familles d'applications d'allure polynomiale \cite{Bi}. Nous proc\'ederons en cinq \'etapes.\\

{\it Normalisations au voisinage d'un param\`etre Misiurewicz.} L'origine \'etant un para-m\`etre  Misiurewicz, notons $\gamma$
la courbe holomorphe fournie par la D\'efinition \ref{DefiMisiu}. On peut, quitte \`a renormaliser, supposer que $\gamma$ est d\'efinie sur $D$.
En rempla\c cant $F$ par $F^{p_0}$, on peut aussi supposer que $p_0=1$.\\
 En outre, une conjugaison par
  $(\la,z)\mapsto (\la,T_{\gamma(\la)}(z))$ o\`u $T_{\gamma(\la)}$ est une famille ad\'equate d'automorphismes de $\Pj^k$ permet de supposer que  $\gamma$ est constant \'egal \`a  $z_1 := \gamma(0)$. D\'esignons alors par $\Omega$ une boule centr\'ee en  $z_1$ et de rayon $r$. Si $\rho$ et $r$ sont pris assez petits, on a la situation suivante :
\begin{itemize}
\item[(i)] $F$ est injective et  uniform\'ement  expansive  sur $D_{\rho}\times \Omega$ : il existe  $K >1$ tel que  $$\forall 
(\la,z)\in D_{\rho}\times \Omega ,\;  d\left(F (\la, z),F(\la, z_1)\right)\ge Kd(z,z_1) ;$$ 
\item[(ii)] $(\la,z_1)\in F^{n_0}(C_F)\Leftrightarrow \la=0.$
\end{itemize}

{\it Construction d'un jeu de contractions satisfaisant les hypoth\`eses de la Section \ref{Sec1}.} Fixons $0<\epsilon \ll 1$ et appliquons la Proposition \ref{PropJC} \`a la famille $F$ et l'ouvert $\Omega$. (Rappelons que $\mu_0(\Omega)>0$ puisque 
$ z_1\in J_0$). Quitte \`a renormaliser $D$, nous obtenons une boule $A$ contenue dans $\Omega$ et charg\'ee par $\mu_0$
ainsi qu'une famille $G_1,\cdots,G_m$ de contractions holomorphes d\'efinies sur $\overline{D}\times \overline{A}$ et satisfaisant les estimations suivantes lorsque $n$ est assez grand :
\begin{itemize}
\item[1)] $G_j(D\times { A}) \subset D \times { A}$,\;$\forall 1\le j\le m$,
\item[2)] $\textrm{dist} \left(G_j(D \times { A}), G_k (D\times { A})\right) >0,\; \forall 1\le j\ne k\le m$,
\item[3)] $e^{-n(\chi_k +3\epsilon)} d(z,z') \le d\left(G_j(\la,z),G_j(\la,z')\right) \le e^{-n(\chi_1-\epsilon)} d(z,z')$ sur $D\times { A}$.
\end{itemize}
Rappelons que $m$ et $n$ sont li\'es par l'in\'egalit\'e $m\ge Cd^{kn}$ et que l'on peut donc, quitte
\`a choisir $n$ assez grand, supposer que 
\begin{eqnarray}\label{aaa}
\frac{\ln m}{n(\chi_k+3\epsilon)} \ge\frac {k\ln d}{\chi_k+3\epsilon} -\epsilon.
\end{eqnarray}
Dimininuer un peu $A$ et renormaliser $D$, permet de les remplacer  par leurs adh\'erences dans les conditions 1),2),3) ci-dessus et d'assurer que les hypoth\`eses de la Section \ref{Sec1} sont satisfaites. \\

{\it Existence d'une courbe holomorphe de cycles $J$-r\'epulsifs dans $D\times A$.}
 Comme $\mu_0(A)>0$, il existe un point $z_0 \in J_0\cap A$ qui est $p_0$-p\'eriodique et r\'epulsif pour $f_0$ (d'apr\`es un th\'eor\`eme de Briend-Duval \cite{BD1}, ces points \'equidistribuent la mesure $\mu_0$). Quitte \`a renormaliser $D$, le th\'eor\`eme des fonctions implicites fournit une courbe
 holomorphe $\sigma : D\to A$ telle que $\sigma(0)=z_0$ et $\sigma(\la)$ est $p_0$-p\'eriodique et r\'epulsif pour $f_\la$
 pour tout $\la \in D$. Il reste \`a \'etudier  l'appartenance de $\sigma(\la)$ \`a $J_\la$. Notons $B(\la)$ une boule centr\'ee en $\sigma(\la)$ et de rayon $r>0$. On peut, quitte \`a renormaliser $D$ et diminuer $r$, supposer que $f_\la^{p_0}$ soit uniform\'ement expansive sur $B(\la)$ pour tout $\la\in D$. Par ailleurs, $\sigma(0)$ appartenant \`a $J_0$ et $\la\mapsto J_\la$ \'etant $s.c.i$ pour la distance de Hausdorff, on peut renormaliser \`a nouveau $D$ de fa\c con \`a ce que $B(\la)\cap J_\la \ne \emptyset$ pour tout $\la \in D$. 
Dans ces conditions, $\sigma(\la)$ est accumul\'e par des pr\'eimages par des it\'er\'ees de $f_\la^{p_0}$ de points de $B(\la)\cap J_\la$   et donc, $J_\la$ \'etant ferm\'e et totalement invariant, $\sigma(\la)\in J_\la$ pour tout $\la\in D$.\\

{\it Mise en place des hypoth\`eses de la Proposition \ref{PropZ}}. En utilisant la seconde \'etape et de la Section \ref{Sec1}, on obtient une famille de graphes holomorphes $\mathcal G$ engendr\'ee par la collection de contractions $G_1,\cdots,G_m$.

Le point $(0,z_1)$ \'etant dans $F^{n_0+n_0'}(C_F)$ pour tout $n_0'\in \N$, notons $Z_0$ la composante irreductible de  $F^{n_0+n_0'}(C_F) \cap (D\times \Omega)$ qui contient $(0,z_1)$. Il est clair que $Z_0$ n'est pas r\'eduite \`a une fibre de $\pi_D$. En tenant compte des normalisations de la premi\`ere \'etape, on voit que si $n_0'$ est pris assez grand alors $\pi_D(Z_0) \Subset D$. Toute  composante irreductible $Z$ de $Z_0\cap(D\times A)$ satisfait   les hypoth\`eses de la Proposition \ref{PropZ}  et donne donc lieu \`a l'estimation
\begin{eqnarray}\label{aaaa}
\textrm{dim}_H \pi_D\left({\mathcal G} \cap Z\right) \ge \frac{\chi_1-\epsilon}{\chi_k +3\epsilon} \left(\frac{\ln m}{n(\chi_k +3\epsilon)}\right) - (2k-2).
\end{eqnarray}

{\it Avalanche de param\`etres Misiurewicz et conclusion}. Nous allons montrer que tous les \'el\'ements de $\pi_D({\mathcal G}\cap Z)$ sont accumul\'es par des param\`etres
Misiurewicz et appartiennent  donc, d'apr\`es le Th\'eor\`emes \ref{thMis}, au lieu de bifurcation. La conclusion r\'esultera alors imm\'ediatement des
estimations (\ref{aaa}) et (\ref{aaaa}) en faisant tendre $\epsilon$ vers $0$.\\
Nous utiliserons la courbe $\sigma$ exhib\'ee \`a la troisi\`eme \'etape. Reprenons les notations de la Section \ref{Sec1}. Par construction, tout $\Gamma_\omega \subset {\mathcal G}$ est  limite, lorsque $p\to +\infty$, de la suite d\'ecroisantes $ T_{\omega_0\cdots\omega_p}= G_{\omega_0\cdots\omega_p}(\overline{D}\times\overline{A})$. De plus cette convergence est uniforme (voir (\ref{a})). En particulier, la suite de graphes $G_{\omega_0\cdots\omega_p}(\Gamma_\sigma)$ converge vers $\Gamma_\omega$. Il s'ensuit que $G_{\omega_0\cdots\omega_p}(\Gamma_\sigma)\cap Z$ converge vers
$\Gamma_\omega \cap Z$ et $\pi_D\left(G_{\omega_0\cdots\omega_p}(\Gamma_\sigma)\cap Z\right)$ vers $\pi_D\left(\Gamma_\omega \cap Z\right)$; d'apr\`es le Lemme \ref{Lem}  ces intersections sont non vides. Or les \'el\'ements de $\pi_D\left(G_{\omega_0\cdots\omega_p}(\Gamma_\sigma)\cap Z\right)$ sont des param\`etres Misiurewicz car la non vacuit\'e de $G_{\omega_0\cdots\omega_p}(\Gamma_\sigma)\cap Z$ signifie que $\Gamma_\sigma$ rencontre une composante de $F^{p+1+n_0+n_0'}(C_F)$ sans y \^{e}tre contenu.\fin

\end{document}